\documentclass[12pt]{article}

\usepackage{amsmath,amssymb,latexsym,amstext,amsopn,amsthm,mathrsfs}

\newtheorem{theorem}{Theorem}
\newtheorem{lemma}{Lemma}

\newtheorem{note}{Remark}

\newcommand{\lqqd}{\hfill{$\Box$}\bigskip}
\newcommand{\CC}{\mathbb{C}}
\newcommand{\NN}{\mathbb{N}}
\newcommand{\RR}{\mathbb{R}}

\newcommand{\dsty}{\displaystyle}
\newcommand{\unifn}{\;\; {\mathop{\mbox{\LARGE
$\rightrightarrows$}}_{n\to\infty}}\;\;}

\begin{document}

\title{On Polar Legendre Polynomials}

\author{H\'{e}ctor Pijeira Cabrera \thanks{Research partially supported by
 Direcci\'{o}n General de Investigaci\'{o}n, Ministerio de Ciencias y Tecnolog\'{\i}a de Espa\~na,
under grant MTM2006-13000-C03-02,  by  Comunidad de Madrid-Universidad Carlos
III de Madrid, under grant CCG06-UC3M/EST-0690 and by Centro de Investigaci\'{o}n Matem\'{a}tica
de Canarias (CIMAC), email: hpijeira@math.uc3m.es}\\
Universidad Carlos III de Madrid, Spain \and Jos\'{e} Y.  Bello Cruz  \thanks{Research  supported
by CNPq-TWAS, email: yunier@impa.br}\\ IMPA, Brazil
 \and Wilfredo Urbina Romero \thanks{Research partially supported by Centro
de
Investigaci\'{o}n Matem\'{a}tica de Canarias (CIMAC), email: wurbina@euler.ciens.ucv.ve} \\
Universidad Central de Venezuela \& University of New Mexico}

\maketitle
\begin{abstract}
We introduce a new class of polynomials $\{P_{n}\}$, that we call polar Legendre polynomials, they 
appear as  solutions of an inverse Gauss problem of equilibrium position  of a field of forces with $n+1$
unit masses. We study algebraic, differential  and asymptotic properties  of this class of
polynomials, that are simultaneously orthogonal with respect to a differential operator and a
discrete-continuous Sobolev type inner product.
\end{abstract}



\section{Introduction}

Let  $\{L_n\}_{\{n \in \NN\}}$ be the monic Legendre polynomials. It is well know
that  $L_n$ satisfies the orthogonality relation
\begin{equation}\label{ortog1}
   \int_{-1}^{1}{L_n(x) x^{k} dx}=0,  \quad   k=0,1,\cdots,n-1\,,
\end{equation}
and  the  differential equation
\begin {equation}\label{LegEcuDif}
-n(n+1)L_n(z)=\left((1-z^2)L_n^{\prime}(z)\right)^{\prime}\,.
\end{equation}

It can be proved, using integration by parts,  that the derivatives of  $\{L_n\}$ satisfy  the following orthogonality condition
\begin{equation}\label{ortog2}
\int_{-1}^{1}{L_{n+1}^{\prime}(x) x^{k}
   (1-x^2)dx}=0,  \quad   k=0,1,\cdots,n-1\,.
\end{equation}

For a fixed complex number $\zeta$, that we are going to call the pole, let us define the $P_{n} =P_{\zeta,n}$ as a monic polynomial, such that
\begin{equation}\label{defpola}
(n+1)\,L_n(z)=\left((z-\zeta)\,P_{n}(z)\right)^{\prime}=P_{n}(z)+(z-\zeta)P_{n}^{\prime}(z)\,,
\end{equation}
$P_{n}$ is called the $n$-th {\em  polar Legendre polynomial}. Obviously,  $P_{n}$ is a monic
polynomial of degree $n$, that by (\ref{ortog1}) and (\ref{defpola}) satisfies the  following
``orthogonality relation''
\begin{equation}\label{ortog3}
\int_{-1}^{1} [P_{n}(x)+(x-\zeta)P_{n}^{\prime}(x)] x^{k} dx =0, \quad \quad k=0,1,\cdots,n-1.
\end{equation}
This type of orthogonality relations generated by differential operators was introduced initially in \cite{ApLoMa02},
where the existence and uniqueness conditions for more general differential expressions were studied in detail.

Observe that the polynomial
\begin{equation}\label{defprimi}
\Pi_{\zeta,n+1}(z)=(z-\zeta)\,P_{n}(z)
\end{equation}
is the primitive of $(n+1)\,L_n(z),$ such that
$\Pi_{\zeta,n+1}(\zeta)=0$;  that is,
\begin {equation}\label{integral}
\Pi_{\zeta,n+1}(z)=(n+1)\int_{\zeta}^{z}{L_n(t) dt}.
\end{equation}

$\Pi_{\zeta,n+1}$ will be  called the {\em primitive Legendre polynomial.} Notice that the properties of
$P_n= P_{\zeta,n}$ and the properties of  $\Pi_{\zeta,n+1}$ are closely related.

It is important to observe that since  the functions that we are considering are entire functions, we can assume
that the definite integrals considered are line integrals defined on  the straight line segment with initial
point in the lower limit of integration and end point in the upper limit of integration.

Now, combining (\ref{LegEcuDif}) and (\ref{defpola}) and integrating  from $\zeta$ to $z$, we have
the fundamental formula
\begin {equation}\label{LegEcuDif1}
n\,(z-\zeta)\,P_{n}(z) =  (1-\zeta^2) L_n^{\prime}(\zeta) -(1-z^2) L_n^{\prime}(z) \,.
\end{equation}

Furthermore, from (\ref{defpola}) it is easy to see that $\Pi_{\zeta,n+1}(z)$ is the $(n+1)$-th
monic orthogonal polynomial with respect to the Sobolev-type inner product (called
``discrete-continuous type'', see \cite{AlPePiRe99})
\begin {equation*}\langle p,q
\rangle = p(\zeta)q(\zeta) +\int_{-1}^{1}\,
p^{\prime}(x)\,q^{\prime}(x)\,dx\,.
\end{equation*}
In \cite{KwLeJu97},  necessary and sufficient conditions under which such
Sobolev-type orthogonal polynomials satisfy a differential equation of spectral type  with polynomial
coefficients are studied.

The localization of critical points of a given class of polynomials has many physical
and geometrical interpretations. Let us consider for instance, a
field of forces given by a system of $n$ masses $m_j,\;(1 \leq j
\leq n )$ at the fixed points $z_j,\;(1 \leq j \leq n )$, that
repels a movable unit mass at $z$ according to the inverse distance law.

Let  $Q_m (z)= \dsty(z-z_1)^{m_1}\cdot(z-z_2)^{m_2} \cdots (z-z_n)^{m_n}$ where $m=m_1+m_2+\cdots+m_n.$ The
logarithmic derivative of $Q_m (z)$ is
\begin{equation}\label{derlog}
 \frac{d(\log (Q_m (z)))}{dz}=\frac{Q_m^{\prime}
(z)}{Q_m (z)}= \frac{m_1}{(z-z_1)}+\frac{m_2}{(z-z_2)}+ \cdots +
\frac{m_n}{(z-z_n)}\,.
\end{equation}
The conjugate of $\dsty \frac{m_j}{(z-z_j)}$ is a vector directed from $z_j$ to $z$, so this vector represents the force at the movable unit mass $z$ due
to a single fixed particle at $z_j$. By (\ref{derlog}) the positions of equilibrium in the field of force coincide with the zeros of $Q_m^{\prime} $, that are not zeros of $Q_m$. In  particular, all multiple 
zeros of $Q_m$ are equilibrium positions. This result is known as Gauss`s theorem (\cite[Theorem 1.2.1, Ch.3]{MiMiRa94}).

Now, we consider the following inverse problem: let
$z^{\prime}_1,\,z^{\prime}_2,\, \cdots,\,z^{\prime}_{n}$ be the
zeros of the orthogonal polynomial $L_n$ which we assume to be  the equilibrium
positions  of a field of forces with $n+1$ unit masses,  one of
which  is given at the point $\zeta$. What is the location of the remaining
masses? Let $P_{n}$ be the monic polynomial whose zeros are the remaining equlibrium positions. By (\ref{defpola}) and (\ref{defprimi})
\begin{equation}\label{RMasses}
\frac{(n+1)\,L_n(z)}{(z-\zeta)P_{n}(z)}=\frac{1}{z-\zeta}+
\frac{P_{n}^{\prime}(z)}{P_{n}(z)}=\frac{\Pi_{\zeta,n+1}^{\prime}(z)}{\Pi_{\zeta,n+1}(z)}\,.
\end{equation}
Then, according to (\ref{derlog}), (\ref{RMasses}) and the above interpretation of the logarithmic derivative,
the location of the remaining unit masses are the zeros of the polynomial $P_{n})$, or
equivalently the poles of (\ref{RMasses}). For this reason we call $P_{n}$ polar polynomial.

The main purpose of this paper is to study some algebraic, differential and analytic properties
of the polar Legendre  polynomials, or equivalently of the primitive Legendre polynomials. The paper is organized as follows. In section 2 we study  the orthogonality relations and  recurrence relations of the polar Legendre polynomials and section 3 is devoted to the study of
the location of zeros and the  asymptotic behavior of zeros and polynomials.

\section{Orthogonality and recurrence relations}

Besides the well known results on orthogonality mentioned in the previous section, we have
the following additional orthogonality relations between  $\{L_n\}$ and $\{P_n\}$,

\begin{theorem}The  polar Legendre
 polynomial $P_{n}$ with pole $\zeta \in \CC$ verifies

\begin {equation}\label{delta1}
\int_{-1}^1 \left[P_{n}(x)+(x-\zeta)P_{n}^{\prime}(x) \right] L_m(x) dx= \left\{
\begin{array}{cc}
       0 & m \neq n, \\   \dsty    (n+1)\|L_n\|^2 & m=n ,\\
\end{array}
\right.
\end{equation} where  $\dsty  \|L_n\|^2=\int_{-1}^1 L^2_n(x) dx. $

Additionally if $n>0$ then
\begin {equation}\label{ortog4}
\int_{-1}^{1} (x-\zeta)P_{n}(x) L_m(x)dx = \left\{
\begin{array}{rl}
 \dsty  \frac{2}{n}
(1-\zeta^2) L_n^{\prime}(\zeta) & m=0, \\
  0 \; & 0 < m < (n-1), \\
\dsty - \frac{n+1}{n}\|L_n\|^2   & m=n-1, \\
0 & m=n ,\\
\dsty  \|L_{n+1}\|^2  & m=n+1, \\
  0 \; & (n+1) < m.
\end{array}
\right.
\end{equation}
\end{theorem}

\proof  Since by (\ref{defpola})
$$\int_{-1}^1 \left[P_{n}(x)+(x-\zeta)P_{n}^{\prime}(x) \right] L_m(x) dx= (n+1)\int_{-1}^1 L_n(x)L_m(x) dx,$$
then (\ref{delta1}) is a direct consequence of (\ref{ortog1}). To prove (\ref{ortog4}), let us  denote
$$I_{n,m} =\int_{-1}^{1} (x-\zeta)P_{n}(x) L_m(x)dx.$$

- The case $m=0$ is a direct consequence of  the fundamental relation (\ref{LegEcuDif1}) and (\ref{ortog2}).

-  If $ 0 < m < (n-1)$  by (\ref{integral}) and Fubini's theorem
\begin{eqnarray}
\nonumber I_{n,m} &=& (n+1)\int_{-1}^{1}\left(\int_{\zeta}^1 L_n(t)dt + \int_{1}^x L_n(t)dt \right)L_m(x)dx  \\
\nonumber   &=& (n+1)\int_{-1}^{1}\left( \int_{1}^x L_n(t)dt \right)L_m(x)dx   \\
\label{SemiOrto}& = & (n+1)\int_{-1}^{1}\left(\int_{t}^{-1}L_m(x) dx \right) L_n(t) dt =0,
\end{eqnarray}
since $\int_{t}^{-1}L_m(x) dx$ is a polynomial (in $t$) of degree $m+1<n$. Notice,  that (\ref{SemiOrto}) is true if we use any other polynomial $Q_m$ of degree  $m<n-1$,
instead of $L_m$.

-If $m=n-1$, observe that as consequence of  (\ref{LegEcuDif}) we have
$$
F_{n}(x)= \int L_{n-1}(x) dx=  \frac{(x^2-1)}{n\,(n-1)}  L_{n-1}^{\prime}(x)= \frac{x^n}{n}+
\ldots\,;
$$
therefore, $F_{n}(-1)=F_{n}(1)=0$. Hence, integrating by parts (\ref{ortog4}), using (\ref{ortog2}) and
(\ref{ortog1}) we get,
$$
I_{n,n-1} = - (n+1)\int_{-1}^{1}L_n(x)F_n(x) dx= -
   \frac{n+1}{n}\|L_n\|^2.$$

-For the case $m=n$,  remember that the Legendre polynomials satisfy  the identity $L_n(x)= (-1)^n L_n(-x)$, $x \in \RR$, and; therefore, the powers of  $L_n$ are either all
odd or all even. By (\ref{integral}) it is clear that  the powers of $(x-\zeta)P_{n}(x)$ have opposite
parity from the ones of $L_n$; consequently  all the powers of $(x-\zeta)P_{n}(x) L_n(x)$ are odd  and since we are integrating in a symmetric interval,
$$ I_{n,n}=  \int_{-1}^{1} (x-\zeta)P_{n}(x) L_n(x)dx=  0,$$

- If  $m=n+1$ by (\ref{ortog1})
$$I_{n,n+1} =\int_{-1}^{1} (x-\zeta)P_{n}(x) L_{n+1}(x)dx= \int_{-1}^{1} x^{n+1} L_{n+1}(x)dx=
\int_{-1}^{1}  L^2_{n+1}(x)dx = \|L_{n+1}\|^2.$$

-Finally, the case $m > (n+1)$ is a direct
consequence  of (\ref{ortog1}).  \lqqd

 As a consequence of these orthogonality relations, let us prove now a recurrence relation for the polar Legendre polynomials,

\begin{theorem} The  polar Legendre
 polynomials $\{P_{n}\}$ with pole $\zeta \in \CC$, satisfy the following  recurrence relation
\begin{equation}\label{recurrencia}
P_{n+1}(z) = z \, P_{n}(z)+ a_{n}P_{n-1}(z)+ b_{n},
\end{equation}
for $n > 1$, where $P_0\equiv 1$, and
$P_1(z)=z+ \zeta$,
\begin {equation}\label{CoefRecurrencia}
a_{n}= \dsty \frac{1-n^2}{n^2} \left( \frac{\|L_{n}\|}{\|L_{n-1}\|}\right)^2 \quad \mbox{and} \quad
b_{n}=\frac{2}{n}(\zeta^2-1)L_n^{\prime}(\zeta).
\end{equation}
\end{theorem}

\proof Let $\{\alpha_{n, 0},\, \alpha_{n, 1}, \, \ldots, \, \alpha_{n, n}\}$ be a set of $n+1$
coefficients such that
\begin{equation}\label{ident}
\dsty (x-\zeta)P_{n}(x)=\sum_{k=0}^{n+1}\alpha_{n, k} P_{k}(x).
\end{equation}
Observe that $\alpha_{n, n+1}=1$ since the polar polynomials $P_n$ are monic. Now in
 order to determine the other coefficients $\alpha_{n, k}, \, k = 0, \cdots, n,$ let us consider

\begin{eqnarray}\label{rerela}
   \nonumber
   (z-\zeta) \left[P_{n}(z)+(n+1)L_n(z) \right] & = &
   (z-\zeta) \left[P_{n}(z)+\left((z-\zeta)P_{n}(z)\right)^{\prime} \right]\\
   \label{RelRec-1} & = &
\sum_{k=0}^{n+1}\alpha_{n, k}  \left[P_{k}(z)+(z-\zeta)P_{k}^{\prime}(z) \right].
\end{eqnarray}

By the orthogonality relation (\ref{ortog4}), we have,
\begin{eqnarray}\label{RecurIzquierdo}
 \sum_{k=0}^{n+1}\alpha_{n, k} \int_{-1}^1   L_m(x)[P_{k}(x) +(x-\zeta)P'_{k}(x)]dx= \alpha_{n, m} \, (m+1) \, \|L_m\|^2,
\end{eqnarray}
for $m=0,\,1,\, \ldots,\, n.$ On the other hand, (\ref{ortog1}) and (\ref{ortog4}) give us
\begin{eqnarray}\label{RecurDerecho}
 \int_{-1}^1
 L_m(x)(x-\zeta)\left[P_{n}(x) +(n+1)L_n(x) \right] dx &=&
 \left\{
\begin{array}{rl}
 \dsty  \frac{2}{n}
(1-\zeta^2) L_n^{\prime}(\zeta) & m=0 \\
  0\quad  & 0 < m < (n-1) \\
\dsty   \frac{n^2-1}{n}\|L_n\|^2   & m=n-1 \\
  - \zeta (n+1) \|L_n\|^2 & m=n.
\end{array}
\right.
\end{eqnarray}

Thus multiplying (\ref{rerela}) by $L_m$, integrating over $[-1,1]$ and  using
(\ref{RecurIzquierdo})--(\ref{RecurDerecho}) we get,

$$\alpha_{n, m}  =\left\{
\begin{array}{cl}
 \dsty  \frac{2}{n}
(1-\zeta^2) L_n^{\prime}(\zeta) & m=0 \\
  0 \; & 0 < m < (n-1) \\
\dsty    \frac{n^2-1}{n^2} \left(\frac{\|L_n\|}{\|L_{n-1}\|}\right)^2   & m=n-1 \\
\dsty   - {\zeta} & m=n.
\end{array}
\right.$$
Replacing these values  in (\ref{ident}) and after cancelations  we get (\ref{recurrencia}). \, \lqqd

\section{Zeros and asymptotics}

Let us now study the distribution of the zeros of the polar Legendre polynomials.  The
next lemma collects several direct consequences of the formulas contained in the introductory
section, in special the fundamental formula (\ref{LegEcuDif1}).

\begin{lemma}\label{ZerosLocat} The  polar Legendre
 polynomials $\{P_{n}\}$ with pole $\zeta \in \CC$, satisfy
\begin{enumerate}
\item If $n$ is odd and $\zeta \in \RR^*=\RR \setminus \{0\}$,
then $x=-\zeta$ is a zero of $P_{n}$.

\item The zeros of the the polar polynomial $P_{n}$ have
multiplicity at most $2$ and the multiple zeros are on
$[-1,1]$.

\item If $\zeta = 1$, $\zeta = -1$ or $L_n^{\prime}(\zeta)=0$, then the
zeros of $P_{n}(x)$ are $-1$ or $1$ and the $(n-1)$ critical points of the Legendre polynomial
$L_n$.

 \item All the zeros of $P_{n}$ are on the
 following  lemniscate
 \begin{equation}\label{lemniscate}
\Lambda_n(\zeta):= \left\{ z\in\CC : \prod_{k=0}^{n}|z-x_{n,k}|=\rho_n(\zeta) \right\},
\end{equation}
where $\dsty \rho_n(\zeta)= \prod_{k=0}^{n}|\zeta-x_{n,k}|$, $x_{n,0}=-1$, $x_{n,n}=1,$ and
$x_{n,1},\, x_{n,2},\, \cdots, \, x_{n,n-1}$ are the $(n-1)$ critical points of the Legendre
polynomial $L_n$.
\end{enumerate}
\end{lemma}

\proof

In order to prove $1.$, using the fact that if $n$ is odd then  all  the powers of  $L_n$ are
odd and (\ref{integral}), one gets $\dsty P_{n}(-\zeta)=0.$

Now, let us suppose  that $w$ is a zero of $P_{n}$ of multiplicity
greater or equal to $3$. Notice that by (\ref{defpola}) a zero of $P_{n}$ with multiplicity greater than 2 is a zero of $L_n$ and also a zero of  $L'_n$, thus $L_n(w)=L^{\prime}_n(w)=0$ which would imply that $w$ is a zero of multiplicity 2 of  $L_n$. This is impossible since the zeros of $L_n$ are all simple, so we have proved $2.$

 Assertion $3.$  is a direct consequence of fundamental formula (\ref{LegEcuDif1}) since  \begin{equation}\label{ceros}
 \left|z_0^2-1\right| \, \left|L_n^{\prime}(z_0)\right|
 = \left|\zeta^2-1\right| \, \left|L_n^{\prime}(\zeta)\right|,
\end{equation}
and then $4.$  follows by considering the factorization of $(z^2-1)L_n^{\prime}(z)$.

 \hfill \lqqd

\begin{note} The following example shows that the zeros of $P_{n}(z)$ do not have to be simple. Let
$\zeta=\frac{2\sqrt{3}}{3}$ (or $\zeta=\,- \frac{2\sqrt{3}}{3}$), hence the corresponding polar Legendre
 polynomial of degree two is
 $$P_{2}(z)=z^2 + \frac{2\sqrt{3}}{3}z+\frac{1}{3}\,  (\mbox{or} \, \,
P_{2}(z)=z^2 - \frac{2\sqrt{3}}{3}z+\frac{1}{3}). $$
Notice that $z= - \frac{\sqrt{3}}{3}$ ( or \, $z=
\frac{\sqrt{3}}{3}$) is a zero of multiplicity two  of $P_{2}(z)$. \\
\end{note}

For the boundedness of the  zeros of the polar Legendre polynomials $\dsty \{P_n\}$ we have the following result,

\begin{lemma} \label{LemaAcotaCeros} Given $\zeta \in \CC$ let us define $\dsty \Delta_{\zeta}= \sup_{x \in [-1,1]}|\zeta-x|$ and
 $\dsty \delta_{\zeta}= \inf_{x \in [-1,1]}|\zeta-x|$, then
\begin{enumerate}
  \item All the zeros of the polar Legendre polynomials $\dsty \{P_n\}$ with pole $\zeta \in \CC$ are contained in $|z| \leq
  \Delta_{\zeta}+1$.
  \item If $\dsty \delta_{\zeta} >1$,
   the zeros of the polar Legendre polynomials $\dsty \{P_n\}$ with pole $\zeta \in \CC$ are simple and contained in the exterior of the ellipse $|z+1|+|z-1|= 2\alpha$, where $ 1< \alpha <\delta_{\zeta}. $
\end{enumerate}
\end{lemma}

\bigskip
\proof

1.  We know by(\ref{lemniscate})  that the zeros of $P_n(z)$ are on the lemniscate
  $\dsty \Lambda_n(\zeta)$. Since $\dsty \rho_n(\zeta) < \Delta_{\zeta}^{n+1}$,  the zeros of $P_n(z)$ are contained in the interior of the lemniscate $\dsty \prod_{k=0}^{n}|z-x_{n,k}|=\Delta_{\zeta}^{n+1}$, where   $x_{n,0}=-1$, $x_{n,n}=1,$ and
$x_{n,1},\, x_{n,2},\, \cdots, \, x_{n,n-1}$ are the $(n-1)$ critical points of the Legendre
polynomial $L_n$, and therefore $|x_{n,k}|\leq 1$. Now, for any $z^*$, such that  $|z^*|> 1+\Delta_{\zeta}$, then
  $$
\prod_{k=0}^{n}|z^*-x_{n,k}| \geq \prod_{k=0}^{n}||z^*|-|x_{n,k}|| > \Delta_{\zeta}^{n+1},
  $$
 and therefore $1.$ is obtained.

2. Let $z$ be such that $|z+1|+|z-1|= 2\alpha$. From the well known arithmetic-geometric mean inequality we get
$$ \prod_{k=0}^{n}|z-x_{n,k}| \leq \left(\frac{1}{n+1} \sum_{k=0}^{n}|z-x_{n,k}| \right)^{{n+1}} < \alpha^{n+1}.$$
If $z$ is a zero of $P_n$, from (\ref{lemniscate}) we get that
$$
\prod_{k=0}^{n}|z-x_{n,k}|=\prod_{k=0}^{n}|\zeta-x_{n,k}| > \delta_{\zeta}^{n+1} > \alpha^{n+1}.
$$
Thus, by these inequalities and $2.$ of lemma  \ref{ZerosLocat} we obtain
 $2.$  of this lemma. \, \lqqd

Finally, let us study the asymptotic behavior of the zeros of the polar Legendre  polynomials,

\begin{theorem}\label{ZerosFuera}
Let $\{P_{n}\} $ be the polar Legendre  polynomials with pole $ \zeta \in \CC
\setminus [-1,1]$, such that $\dsty \delta_{\zeta} >1$. Then the set of accumulation points of
zeros of $\{P_{n}\}$ are on the ellipse
\begin{equation}\label{ZerosAcum}
    \Lambda(\zeta):= \left\{ z \in \CC: z= \frac{\rho^2(\zeta)+1}{2\rho(\zeta)} \cos \theta
    + i \,\frac{\rho^2(\zeta)-1}{2\rho(\zeta)} \sin \theta , \quad 0 \leq \theta  < 2 \pi \right\},
   \end{equation} where $\rho(\zeta):=|\zeta+\sqrt{\zeta^2-1}|$ and the branch of the square
root is chosen so that $\dsty |z+\sqrt{z^2-1}|>1$ for $z \in \CC \setminus [-1,1]$.
\end{theorem}

\proof By (\ref{ceros}) we have that the zeros of the $n$th  polar Legendre polynomial satisfies
the equation
\begin{equation}\label{RaizCeros}
 \left|\frac{z^2-1}{n}\right|^{\frac{1}{n}} \, \left|L_n^{\prime}(z)\right|^{\frac{1}{n}}
 = \left|\frac{\zeta^2-1}{n}\right|^{\frac{1}{n}}
 \, \left|L_n^{\prime}(\zeta)\right|^{\frac{1}{n}}\,.
\end{equation}

On  other hand,  from the asymptotic properties of the Legendre polynomials  it is well known that
\begin{equation}\label{Limelegendre}
    \dsty \lim_{n  \rightarrow
\infty}\left|L^{\prime}_n(z)\right|^{\frac{1}{n}}\,= \frac{|z+\sqrt{z^2-1}|}{2}\,,
\end{equation}
uniformly on compact subset of $\CC \setminus [-1,1]$.  Taking limit as $n \rightarrow \infty$, from 2. of lemma \ref{LemaAcotaCeros}, it is possible to use (\ref{Limelegendre}) on  both sides of (\ref{RaizCeros}) and we have
that set of accumulation points of the zeros of the sequence of polynomials $\{P_n\}$ are contained
in the curve $$\Lambda(\zeta)=\left\{ z \in \CC: |z+\sqrt{z^2-1}|=\rho(\zeta)\right\}\,. $$
Hence
\begin{eqnarray*}
z+\sqrt{z^2-1}&=& \rho(\zeta)\, e^{i\theta}, \quad 0 \leq \theta  < 2 \pi,\\
z-\sqrt{z^2-1}&=& \rho(\zeta)^{-1}\,e^{-i\theta}\,,\\
2z & = & \rho(\zeta) e^{i\theta}+ \rho(\zeta)^{-1}e^{-i\theta}.
\end{eqnarray*} \lqqd

Finally, let us study the relative asymptotic of the polar Legendre polynomials $\{ P_n \}$ with respect to 
the Legendre polynomials $\{ L_n \}$ and their derivatives $\{ L'_n \}$

\begin{theorem}\label{RelativeAsymptotic}
Let $\zeta$ be a fixed complex number and  $L_n$, $P_n$ y $\Delta_{\zeta}$ be as above. Then
\begin{eqnarray}
 \label{AsintComp}  1.  & \dsty \frac{n\,P_{n}(z)}{L_n^{\prime}(z)} & \unifn  \frac{z^2-1}{z-\zeta}, 
\mbox{uniformly on compact subsets  of the set }\\
\nonumber   & &  \{z \in \CC : |z| > \Delta_{\zeta}+1\}.\\
\label{AsintComp1} 2.  &\dsty \frac{P_{n}(z)}{L_n(z)}  &  \unifn  \frac{\sqrt{z^2-1}}{z-\zeta}, 
\mbox{uniformly on compact subsets  of the set }\\
\nonumber   & &  \{z \in \CC : |z|
> \Delta_{\zeta}+1\}.
\end{eqnarray}
\end{theorem}

\proof Let $K$ be a compact subset of $\{z \in \CC : |z| > \Delta_{\zeta}+1\}$. By the fundamental
formula (\ref{LegEcuDif1}) we have
$$
\frac{n\,P_{n}(z)}{L_n^{\prime}(z)} =  \frac{1-\zeta^2}{z-\zeta}
\frac{L_n^{\prime}(\zeta)}{L_n^{\prime}(z)} +\frac{z^2-1}{z-\zeta} \,.$$
Hence, in order to prove (\ref{AsintComp}) it is sufficient to show that
\begin{equation}\label{SufAsint}
\frac{L_n^{\prime}(\zeta)}{L_n^{\prime}(z)} \unifn 0
\end{equation}
uniformly on compact subset $K \subset \{z \in \CC : |z| > \Delta_{\zeta}+1\}$. Let
$x_{n,1}^{\prime},\,x_{n,2}^{\prime}, \dots, x_{n,n-1}^{\prime}$ be the  $n-1$ zeros of
$L_n^{\prime}(z)$ and$\dsty d_{\zeta,K}= \inf_{{z \in K}\atop{|w|=\Delta_{\zeta}+1}}|z-w|$, hence
$$ \left| \frac{L_n^{\prime}(\zeta)}{L_n^{\prime}(z)} \right|=
\frac{\dsty \prod_{k=1}^{n-1}\left|\zeta-x_{n,k}^{\prime}\right|}{\dsty
\prod_{k=1}^{n-1}\left|z-x_{n,k}^{\prime}\right|}
 < \left( \frac{\Delta_{\zeta}}{d_{\zeta,K}+\Delta_{\zeta}}\right)^{n-1} < 1, \quad z \in K.
$$ which is equivalent to the uniform convergence of (\ref{SufAsint}) to zero on  a compact subset
$K$ of $ \{z \in \CC : |z| > \Delta_{\zeta}+1\}$.

The statement (\ref{AsintComp1}) is a direct consequence of (\ref{AsintComp}) and the well known
asymptotic behavior of the Legendre polynomials (see \cite[corollary 1.6]{Van87}) $$
\frac{L_n^{\prime}(z)}{n\,L_{n}(z)}\unifn \frac{1}{\sqrt{z^2-1}}, \; \mbox{uniformly on compact
subsets of }\; \CC \setminus [-1,1].
$$ \lqqd

\end{document}